\documentclass{article}%
\usepackage{amsmath}
\usepackage{amsfonts}
\usepackage{amssymb}
\usepackage{graphicx}
\usepackage{booktabs}
\usepackage{subfigure}
\usepackage{multirow}
\usepackage{array}
\usepackage{hhline}
\usepackage{setspace}
\usepackage{hyperref}
\usepackage{latexsym}
\usepackage{psfrag}
\usepackage{subfigure}
\hypersetup{colorlinks=true, linkcolor=blue, citecolor=red, urlcolor=red}
\newtheorem{theorem}{Theorem}

\newenvironment{proof}[1][Proof]{\noindent\textbf{#1.} }{\ \rule{0.5em}{0.5em}}

\begin{document}
\doublespacing
\title{ \textbf{Measuring Bayesian Robustness Using R\'enyi Divergence}}   

\author{Luai Al-Labadi \thanks{{\em Address for correspondence}: Department of Mathematical \& Computational Sciences, University of Toronto Mississauga, Ontario L5L 1C6, Canada. E-mail:  luai.allabadi@utoronto.ca.}, Forough Fazeli Asl \thanks{Department of Mathematical Sciences, Isfahan University of Technology,
Isfahan 84156-83111, Iran. E-mail: forough.fazeli@math.iut.ac.ir.}, \& Ce Wang\thanks {Department of Mechanical and Industrial Engineering, University of Toronto, Canada. E-mail: ce.wang@mail.utoronto.ca.}}

\date{\vspace{-5ex}}
\maketitle
\pagestyle {myheadings} \markboth {} {Measuring Robustness Using R\'enyi Divergence }
\begin{abstract}
This paper deals with measuring the Bayesian robustness of classes of contaminated priors. Two different classes of priors in the neighborhood
of the elicited prior are considered. The first one is the well-known $\epsilon$-contaminated class, while the second one is the geometric mixing class.
The proposed measure of robustness is based on computing the curvature of R\'enyi divergence between posterior distributions.  Examples are used to illustrate the results by using simulated and real data sets.
\par
 \vspace{9pt} \noindent\textsc{Keywords:}  Bayesian Robustness, $\epsilon$-contamination, Geometric contamination,  R\'enyi divergence.

 \vspace{9pt}

\noindent { \textbf{MSC 2000}} 62F15, 62F35, 62G35.
\end{abstract}

\section{Introduction}  Bayesian inferences require the specification of a prior, which contains a priori knowledge about the parameter(s). If the selected prior, for instance, is flawed, this may yield erroneous inferences.

The goal of this paper is to measure the sensitivity of inferences to a chosen prior (known as \emph{robustness}). Since, in most cases, it becomes very challenging to come up with only a sole prior distribution, we consider a class, $\Gamma$, of all possible priors over the parameter space. To construct $\Gamma$, a preliminary prior $\pi_0$ is elicited. Then robustness for all priors $\pi$ in a neighborhood of $\pi_0$ is intended. A common accepted way to construct neighborhoods around $\pi_0$ is through contamination.  Specifically, we will consider two different classes of contaminated or mixture of priors, which are given by
\begin{equation}
\label{contaminated} \Gamma_a=\left\{\pi(\theta): \pi(\theta)=(1-\epsilon)\pi_0(\theta)+\epsilon q(\theta), q \in Q\right\}
\end{equation}
and
\begin{equation}
\label{geometric} \Gamma_g=\left\{\pi(\theta): \pi(\theta)=c(\epsilon) \pi_0^{1-\epsilon}(\theta)q^{\epsilon}(\theta), q \in Q\right\},
\end{equation}
where $\pi_0$ is the elicited prior, $Q$ is a class of distributions, $c(\epsilon)$ is normalizing constant and $0 \le  \epsilon \le 1$ is a small given number denoting the amount of contamination. For other possible classes of priors, see for instance, De Robertis and Hartigan (1981) and Das Gupta and Studden (1988a, 1988b).

The class (\ref{contaminated}) is known as the $\epsilon$-contaminated class of priors. Many papers about the class (\ref{contaminated}) are found in the literature. For instance, Berger (1984, 1990), Berger and Berliner (1986), and Sivaganesan and Berger (1989) used various choices of $Q$.  Wasserman (1989) used (\ref{contaminated}) to study robustness of likelihood regions. Dey and Birmiwal (1994) studied robustness based on the curvature. Al-Labadi and Evans (2017) studied robustness of relative belief ratios (Evans, 2015)  under class (\ref{contaminated}).

On the other hand, the class (\ref{geometric}) will be referred as geometric contamination or mixture class. This class was first studied, in the context of Bayesian Robustness,  by Gelfand and Dey (1991), where the posterior robustness was measured using Kullback-Leibler divergence.  Dey and Birmiwal (1994) generalized the results of Gelfand and Dey (1991)  under (\ref{contaminated}) and  (\ref{geometric}) by using the $\phi$-divergence defined by $d(\pi(\theta|x), \pi_{0}(\theta|x))=\int \pi_{0}(\theta|x) \phi (\pi(\theta|x)/\pi_{0}(\theta|x))d \theta$ for a smooth convex function $\phi$. For example, $\phi(x)=x\ln x$ gives Kullbak-Leibler divergence.

In this paper, we extend the results of Gelfand and Dey (1991) and Dey and Birmiwal (1994) by applying R\'enyi divergence on both classes (\ref{contaminated}) and (\ref{geometric}). This will give local sensitivity analysis on the effect of small perturbation to the prior. R\'enyi entropy, developed by Hungarian mathematician Alfr\'ed  R\'enyi in 1961, generalizes the Shannon entropy and includes other entropy measures as special cases. It  finds applications, for instance, in statistics (Kanaya and Han, 1995), pattern
recognition (Jenssen, Hild, Erdogmus, Principe and Eltoft, 2003), economics (Bentes, Menezes and Mendes, 2008) and biomedicine (Lake, 2006).

An outline of this paper is as follows. In Section 2, we give definitions, notations and some properties of R\'enyi
divergence. In Section 3, we develop curvature formulas for measuring robustness
based on R\'enyi divergence. In Section 4, three examples are studied to illustrate the results numerically.
Section 5 ends with a brief summary of the results.

\section{Definitions and  Notations}  Suppose we have a statistical model that is given by the density function $f_\theta(x)$ (with respect to some measure), where $\theta$ is an unknown parameter that belongs to the parameter space $\Theta$. Let $\pi(\theta)$ be the prior distribution of $\theta$. After observing the data $x$, by Bayes' theorem, the posterior distribution of $\theta$  is given by the density
\begin{equation*}
  \pi(\theta|x) = \frac{f_\theta(x)\pi(\theta)}{m(x|\pi)},
\end{equation*}
where
\begin{equation*}
m(x|\pi) = \int f_\theta(x)\pi(\theta) d\theta
\end{equation*}
is the prior predictive density  of the data.

To measure the divergence between two posterior distributions, we consider  R\'enyi divergence (R\'enyi, 1961). R\'enyi divergence of order $a$  between two posterior densities $\pi(\theta|x)$ and $\pi_0(\theta|x)$ is defined as:
\begin{eqnarray*}
 d=d(\pi(\theta|x),\pi_{0}(\theta|x))&=&\frac{1}{a-1}\ln\left(\int{{{\left(\pi(\theta|x)\right)^a}{\left(\pi_{0}(\theta|x)\right)^{1-a}}d\theta}}\right)\nonumber\\
                               &=&\frac{1}{a-1}\ln \left(E_{\pi_{0}(\theta|x)}\left[\left(\frac{\pi(\theta|x)}{\pi_{0}(\theta|x)}\right)^a\right]\right),\label{renyi}
\end{eqnarray*}
where $a>0$ and $E_{\pi_{0}(\theta|x)}$ denotes the expectation with respect to the density $\pi_0(\theta|x)$. It is known that $d(\pi(\theta|x),\pi_{0}(\theta|x))\ge 0$ for all $\pi(\theta|x),\pi_{0}(\theta|x), a>0$ and   $d(\pi(\theta|x),\pi_{0}(\theta|x))= 0$ if and only if $\pi(\theta|x)=\pi_{0}(\theta|x)$. Note that, the case  $a = 1$ is defined by letting $a \to 1$. This leads to the Kullbak-Leibler divergence. For further properties of  R\'enyi divergences consult, for example, Li and Turner (2016).

Following  the idea of McCulloch (1989) and Dey and Birmiwal (1994) for calibrating, respectively,  the Kullback-Leibler  divergence and the $\phi$ divergence, it is also possible to calibrate R\'enyi divergence as follows. Consider  a biased coin where $X=1$ (heads) occurs with probability $p$.  Then R\'enyi divergence between an unbiased and a biased coin is
 \begin{equation*}
    d(f_0,f_1)=\frac{1}{a-1}\ln\left[2^{a-1}\left(p^a+(1-p)^{a}\right)\right],
\end{equation*}
where, for $x=0,1$, $f_0(x)=0.5$ and $f_1(x)=p^x(1-p)^{1-x}$. Now, setting $d(f_0,f_1)=d_0$ 
 gives
\begin{equation}
   2^{1-a}e^{(a-1)d_0}=p^a+(1-p)^{a}. \label{interpret1}
\end{equation}
 Then the number $p$ is the calibration of $d$. In general, equation (\ref{interpret1}) needs to be solved numerically for $p$. Note that, for the case $a=1$ (i.e. the Kullback-Leibler divergence) one may use the following explicit formula for $p$ due to McCulloch (1989):
 \begin{equation}
  p=0.5+0.5\left(1-e^{-2d_0}\right)^{1/2}. \label{interpret2}
\end{equation}
Values of $p$ close to 1 indicate that $f_0$ and $f_1$ are quite different, while values of $p$ close to 0.5 implies that they are similar. It is restricted that $p$ is chosen so that it is between 0.5 and 1 there is a one-to-one correspondence between $p$ and $d_0$. 
%
%
%

A motivating key fact about  R\'enyi divergence follows from its Taylor expansion. Let
$$f(\epsilon)=d(\pi(\theta|x),\pi_{0}(\theta|x))=
\frac{1}{a-1}\ln\left(\int{{{\left(\pi(\theta|x)\right)^a}{\left(\pi_{0}(\theta|x)\right)^{1-a}}d\theta}}\right),$$ where $\pi(\theta|x)$ is the posterior distribution of $\theta$ given the data $x$ under the prior $\pi$ defined in (\ref{contaminated}) and (\ref{geometric}). Assuming differentiability with respect to $\epsilon$, the Taylor expansion of $f(\epsilon)$ about  $\epsilon=0$ is given by
 \begin{equation*}
  f(\epsilon)=f(0)+\epsilon \frac{\partial f(\epsilon)}{\partial \epsilon} \bigg|_{\epsilon=0}+\frac{\epsilon^2}{2} \frac{\partial^2 f(\epsilon)}{\partial\epsilon^2}\bigg|_{\epsilon=0}+ \cdots.
\end{equation*}
Clearly, $f(0)=0$.  If integration and differentiation are interchangeable, we have
\begin{eqnarray*}
 \frac{\partial f(\epsilon)}{\partial \epsilon}&=&\frac {a}{1-a}\frac{\int \left(\pi_0(\theta|x)\right)^{1-a} \left(\pi(\theta|x)\right)^{a-1} \frac{\partial \pi(\theta|x)}{\partial \epsilon} d\theta}
 {\int \left(\pi_0(\theta|x)\right)^{1-a} \left(\pi(\theta|x)\right)^{a} d\theta}.
\end{eqnarray*}
Hence,
\begin{eqnarray*}
\nonumber \frac{\partial f(\epsilon)}{\partial \epsilon}\bigg|_{\epsilon=0}&=&\frac {a}{1-a} \int  \frac{\partial\pi (\theta|x)}{\partial \epsilon}  d \theta\\
                                                                  &=&\frac {a}{1-a} \frac{\partial}{\partial \epsilon} \left(\int \pi (\theta|x) d \theta\right)=\frac {a}{1-a} \frac{\partial}{\partial \epsilon} (1)=0. \label{eq1}
\end{eqnarray*}
On the other hand,

\begin{eqnarray*}
 \frac{\partial^2 f(\epsilon)}{\partial \epsilon^2}&=& \frac{\partial }{\partial \epsilon} \left(\frac {a}{1-a}\frac{\int \left(\pi_0(\theta|x)\right)^{1-a} \left(\pi(\theta|x)\right)^{a-1} \frac{\partial \pi(\theta|x)}{\partial \epsilon} d\theta}
 {\int \left(\pi_0(\theta|x)\right)^{1-a} \left(\pi(\theta|x)\right)^{a} d\theta}\right),
\end{eqnarray*}
which, at $\epsilon=0$, reduces to

\begin{eqnarray*}
 \frac{\partial^2 f(\epsilon)}{\partial\epsilon^2}\bigg|_{\epsilon=0}&=& -a \int \frac{\left(\frac{\partial\pi (\theta|x)}{\partial \epsilon}\right)^2}{\pi(\theta|x)} d \theta \bigg|_{\epsilon=0}\\
                                                                      &=& -a \int \left(\frac{\frac{\partial\pi (\theta|x)}{\partial \epsilon}}{\pi(\theta|x)}\right)^2 \pi(\theta|x)  d \theta \bigg|_{\epsilon=0}\\
                                                                      &=& -a E_{\pi(\theta|x)}\left[\left(\frac{\partial \ln \pi (\theta|x)}{\partial \epsilon}\right)^2 \right]\bigg|_{\epsilon=0}\\
                                                                      &=& -a I_{\pi(\theta|x)}(\epsilon)\bigg|_{\epsilon=0}.
\end{eqnarray*}
Here $I_{\pi(\theta|x)}(\epsilon)=E_{\pi(\theta|x)}\left[\left(\frac{\partial \ln \pi (\theta|x)}{\partial \epsilon}\right)^2 \right]\bigg|_{\epsilon=0}$ is the Fisher information function for $\pi(\theta|x)$ (Lehmann and  Casella, 1998). Thus, for $\epsilon \approx 0$, we have
 \begin{equation}
d(\pi(\theta|x),\pi_{0}(\theta|x)) \approx -\frac{a\epsilon ^2} {2} I_{\pi(\theta|x)}(\epsilon). \label{fisher}
\end{equation}
Note that, $\partial^2 f(\epsilon)/\partial \epsilon^2\bigg|_{\epsilon=0}=\partial^2 d/\partial \epsilon^2 \bigg|_{\epsilon=0}$ is known as the local \emph{curvature}  at $\epsilon=0$  of  R\'enyi divergence. Formula (\ref{fisher}) justifies the use of the curvature to measure the Bayesian robustness of the two classes of priors $\Gamma_a$ and $\Gamma_g$ as defined in (\ref{contaminated}) and (\ref{geometric}), respectively. Also this formula provide a direct relationship between Fisher's information and  the curvature   of  R\'enyi divergence.

\section{Measuring Robustness Using R\'enyi Divergence} In this section, we explicitly obtain the local curvature  at $\epsilon=0$  of  R\'enyi divergence (i.e. ${\partial^2 d}/{\partial \epsilon^2}\bigg|_{\epsilon=0}$), to measure the Bayesian robustness of the two classes of priors $\Gamma_a$ and $\Gamma_g$ as defined in (\ref{contaminated}) and (\ref{geometric}), respectively. The resulting quantities are presumably much easier to estimate than working directly with R\'enyi divergence.

\begin{theorem} \label{theorem1} For the $\epsilon$-contaminated class defined in (\ref{contaminated}), the local curvature of R\'enyi divergence at $\epsilon=0$ is
\begin{equation*}
  C_{a}^{\Gamma_{a}} =\frac{\partial^2 d}{\partial \epsilon^2}\bigg|_{\epsilon=0} =aVar_{\pi_0(\theta|x)}
  \bigg[
  \frac{q(\theta)}{\pi_0(\theta)}
  \bigg],
\end{equation*}
where $Var_{\pi_0(\theta|x)}$ denotes the variance with respect to $\pi_0(\theta|x)$.
\end{theorem}

\proof Under the prior $\pi$ defined in (\ref{contaminated}),  the marginal $m(\theta|x)$ and the posterior distribution $\pi(\theta|x)$ can be written as
\begin{eqnarray*}
 m(x|\pi)=(1-\epsilon) m(x|\pi_0)+\epsilon m(x|q)
   \end{eqnarray*}
and
\begin{eqnarray}
\pi(\theta|x)&=&\frac{f_{\theta}(x) \pi(\theta)}{m(x|\pi)}\nonumber\\
&=& \frac{f_{\theta}(x) \left((1-\epsilon) \pi_{0}(\theta)+\epsilon q(\theta)\right)}{m(x|\pi)}\nonumber\\
              &=& \lambda(x) \pi_0(\theta|x)+(1-\lambda(x)) q(\theta|x),\label{posterior1}
   \end{eqnarray}
where \begin{equation*}
\lambda(x)=(1-\epsilon)\frac{m(x|\pi_0)}{m(x|\pi)}.
\end{equation*}
Define
\begin{eqnarray*}
    f(\epsilon)&=&d\left(\pi(\theta|x),\pi_{0}(\theta|x)\right)\\
    &=&\frac{1}{a-1}\ln\left[\int \left(\pi(\theta|x)\right)^a\left(\pi_0(\theta|x)\right)^{1-a} d\theta\right]=\frac{1}{a-1}\ln\left[\int \gamma d\theta\right], 
\end{eqnarray*}
where
\begin{equation*}
\gamma = \left(\pi(\theta|x)\right)^a\left(\pi_0(\theta|x)\right)^{1-a}=\left(\lambda(x) \pi_0(\theta|x)+(1-\lambda(x)) q(\theta|x)\right)^a\left(\pi_0(\theta|x)\right)^{1-a}.
\end{equation*}
Clearly,
\begin{equation}
\gamma\bigg|_{\epsilon=0} = \pi_0(\theta|x)  \ \ \mbox{and} \ \ \int\gamma\bigg|_{\epsilon=0}d\theta = 1.  \label{gamma}
\end{equation}
We have
\begin{equation*}
\frac{\partial\gamma}{\partial\epsilon}=a\frac{m(x|q) m(x|\pi_0)\left(q(\theta|x)-\pi_0(\theta|x)\right)}
{\left[\epsilon q(\theta|x) m(x|q)+(1-\epsilon)m(x|\pi_0)\pi_0(\theta|x)\right]\left[(1-\epsilon) m(x|\pi_0)+\epsilon m(x|q)\right]}
\end{equation*}
and
\begin{equation*}
\frac{\partial\gamma}{\partial\epsilon}\bigg|_{\epsilon = 0}=a\frac{m(x|q)\left(q(\theta|x)-\pi_0(\theta|x)\right)}{m(x|\pi_0)}.
\end{equation*}
Thus,
\begin{equation}
\int\frac{\partial\gamma}{\partial\epsilon}d\theta\bigg|_{\epsilon = 0}=0. \label{eq6}
\end{equation}
Now,
\begin{align*}
\frac{\partial^2d}{\partial \epsilon^2}=\frac{\partial}{\partial \epsilon}\left(\frac{1}{a-1} \frac{\int\frac{\partial\gamma}{\partial\epsilon}d\theta}{\int\gamma d\theta}\right)&=\frac{1}{a-1}\frac{[\int\gamma d\theta][\int\frac{\partial^2\gamma}{\partial\epsilon^2}d\theta]-[\int\frac{\partial\gamma}{\partial\epsilon}d\theta]^2}{[\int\gamma d\theta]^2}.
\end{align*}
By (\ref{gamma}) and (\ref{eq6}),
\begin{align*}
\frac{\partial^2d}{\partial \epsilon^2}\bigg|_{\epsilon = 0}&=\frac{1}{a-1}\int\frac{\partial^2\gamma}{\partial\epsilon^2}\bigg|_{\epsilon=0}d\theta.
\end{align*}
We have
\begin{align}
\begin{split}
\frac{\partial^2\gamma}{\partial\epsilon^2}\mid_{\epsilon=0}
=&\bigg(\frac{\pi_0(\theta|x)m(x|\pi_0)-q(\theta|x)m(x|q)}{\pi_0(\theta|x)m(x|\pi_0)}+\frac{m(x|\pi_0)-m(x|q)}{m(x|\pi_0)}+\\
&\frac{a\frac{m(x|q)}{m(x|\pi_0)}\left(q(\theta|x)-\pi_0(\theta|x)\right)}{\pi_0(\theta|x)}\bigg)\times\\
&\ a\frac{m(x|q)}{m(x|\pi_0)}\left(q(\theta|x)-\pi_0(\theta|x)\right). \label{2ndderivativ}
\end{split}
\end{align}
Since
\begin{align}
\nonumber  \frac{m(x|q)}{m(x|\pi_0)}&=\frac{\int f_{\theta}(x) q(\theta)d\theta}{m(x|\pi_0)}=\frac{\int f_{\theta}(x) \pi_{0}(\theta)\frac{q(\theta)}{\pi_{0}(\theta)}d\theta}{m(x|\pi_0)} \\
\nonumber  &=\int \pi_0(\theta|x)\frac{q(\theta)}{\pi_{0}(\theta)}d\theta\\
&=E_{\pi_0(\theta|x)}\left[\frac{q(\theta)}{\pi_0(\theta)}\right], \label{expectation}
\end{align}
from (\ref{2ndderivativ}), we get
\begin{align*}
\begin{split}
\frac{\partial^2\gamma}{\partial\epsilon^2}\bigg|_{\epsilon=0}
=&a\left(2-E_{\pi_0(\theta|x)}\left[\frac{q(\theta)}{\pi_0(\theta)}\right]\right)E_{{\pi_0}(\theta|x)}\left[\frac{q(\theta)}{\pi_0(\theta)}\right]
\left(q(\theta|x)-\pi_0(\theta|x)\right)\\
&-a\left(E_{\pi_0(\theta|x)}\right)^2\left[\frac{q(\theta)}{\pi_0(\theta)}\right]\left(\frac{q(\theta|x)}{\pi_0(\theta|x)}\right) \left(q(\theta|x)-\pi_0(\theta|x)\right)\\
&+{a^2}\left(E_{\pi_0(\theta|x)}\right)^2\left[\frac{q(\theta)}{\pi_0(\theta)}\right]\frac{\left(q(\theta|x)-\pi_0(\theta|x)\right)^2}{\pi_0(\theta|x)}.
\end{split}
\end{align*}
Therefore,
\begin{eqnarray}
\nonumber \frac{{\partial^2}d}{\partial\epsilon^2}\bigg|_{\epsilon=0}=a\bigg(\left(E_{\pi_0(\theta|x)}\left[\frac{q(\theta)}{\pi_0(\theta)}\right]\right)^2
E_{\pi_0(\theta|x)}\left[\left(\frac{q(\theta|x)}{\pi_0(\theta|x)}\right)^2\right]\\
-\left(E_{\pi_0(\theta|x)}\left[\frac{q(\theta)}{\pi_0(\theta)}\right]\right)^2\bigg). \label{eq9}
\end{eqnarray}
Note that,
\begin{equation*}
\left(\frac{q(\theta|x)}{\pi_0(\theta|x)}\right)^2=\left(\frac{q(\theta) f_{\theta}(x)/m(x|q)}{\pi(\theta)f_{\theta}(x)/m(x|\pi_0)}\right)^2
=\left(\frac{q(\theta)}{\pi(\theta)}\right)^2 \left(\frac{m(x|\pi_0)}{m(x|q)}\right)^2
\end{equation*}
Hence, by (\ref{expectation}),
\begin{equation}
E_{\pi_0(\theta|x)}\left[\left(\frac{q(\theta|x)}{\pi_0(\theta|x)}\right)^2\right]=
E_{\pi_0(\theta|x)}\left[\left(\frac{q(\theta)}{\pi_0(\theta)}\right)^2\right]\frac{1}
{\left(E_{\pi_0(\theta|x)}\left[\frac{q(\theta)}{\pi_0(\theta)}\right]\right)^2}.  \label{eq10}
\end{equation}
Thus, by (\ref{eq9}) and (\ref{eq10}),
\begin{align*}
\frac{{\partial^2}d}{\partial{\epsilon}^2}\bigg|_{\epsilon=0}&=
a\left(E_{\pi_0(\theta|x)}\left[\left(\frac{q(\theta)}{\pi_0(\theta)}\right)^2\right]-
\left(E_{\pi_0(\theta|x)}\left[\frac{q(\theta)}{\pi_0(\theta)}\right]\right)^2\right)\\
&=aVar_{\pi_0(\theta|x)}\left[\frac{q(\theta)}{\pi_0(\theta)}\right].
\end{align*}
\endproof

\begin{theorem} \label{theorem2}
 For the geometric contaminated class defined in (\ref{geometric}), the local curvature of R\'enyi divergence at $\epsilon$ = 0 is\\
\begin{equation*}
C_{a}^{\Gamma_{g}}=\frac{\partial^2 d}{\partial \epsilon^2}\bigg|_{\epsilon=0} =aVar_{{\pi_0}(\theta|x)}\left[\ln\left(\frac{q(\theta)}{\pi_0(\theta)}\right)\right],
\end{equation*}
$Var_{\pi_0(\theta|x)}$ denotes the variance with respect to $\pi_0(\theta|x)$.
\end{theorem}
\proof
Define
\begin{equation*}
\gamma = \left(\pi(\theta|x)\right)^a\left(\pi_0(\theta|x)\right)^{1-a}.
\end{equation*}
Thus,
\begin{eqnarray*}
 d&=&\frac{1}{a-1}\ln\left(\int \gamma d \theta\right).
\end{eqnarray*}
We have
\begin{eqnarray*}
   \frac {\partial d}{\partial \epsilon}&=&\frac{1}{a-1} \times \frac{\int{\frac {\partial \gamma }{\partial \epsilon}d \theta}} {\int \gamma d \theta}
\end{eqnarray*}
and
\begin{eqnarray}
   \frac {\partial^2 d}{\partial \epsilon^2}&=&\frac{1}{a-1} \times \frac{\int \gamma d \theta \int{\frac {\partial^2 \gamma }{\partial \epsilon^2}d \theta}-\left(\int \frac {\partial \gamma}{\partial \epsilon}  d \theta\right)^2} {\left(\int \gamma d \theta\right)^2}. \label{2nd}
\end{eqnarray}
Since $\gamma\bigg|_{\epsilon=0}=\pi_0(\theta|x)$,
\begin{eqnarray*}
\frac{{\partial^2}d}{\partial{\epsilon^2}}\bigg|_{\epsilon=0}&=&\int\frac{\partial^2\gamma}{\partial\epsilon^2}d\theta\bigg|_{\epsilon=0}
-\left(\int\frac{\partial\gamma}{\partial\epsilon }d\theta\right)^2\bigg|_{\epsilon=0}.
\end{eqnarray*}
For the geometric class defined in (\ref{geometric}),
\begin{eqnarray}\label{posterior2}
\pi(\theta|x)=\frac{f_\theta(x) \pi(\theta)}{m(x|\pi)}=\frac{f_\theta(x) c(\epsilon) (\pi_0(\theta))^{1-\epsilon} (q(\theta))^{\epsilon}}{m(x|\pi)}
 \ \ \mbox{and} \ \ \pi_0(\theta|x)=\frac{f_\theta(x) \pi_0(\theta)}{m(x|\pi_0)}.
\end{eqnarray}
Thus,
\begin{eqnarray*}
\gamma&=&\frac{f_\theta(x)(c(\epsilon))^a (\pi_0(\theta))^{1-a\epsilon}(q(\theta))^{a\epsilon}}{(m(x|\pi))^{a}(m(x|\pi_0))^{1-a}}.
\end{eqnarray*}
Therefore,
\begin{eqnarray*}
\ln\left(\gamma\right)&=&a\ln \left(\frac{c(\epsilon)}{m(x|\pi)}\right)- a\epsilon\ln \left(\frac{\pi_0(\theta)}{q(\theta)}\right)+\ln \frac{f_{\theta}(x) \pi_0(\theta)}{ (m(x|\pi_0))^{1-a}}.
\end{eqnarray*}
We have
\begin{eqnarray}
\frac{\partial \gamma}{\partial \epsilon}=\gamma \frac{\partial \ln \gamma}{\partial \epsilon}=a \gamma\left(\frac{{\partial}}{\partial{\epsilon}}\ln \left(\frac{c(\epsilon)}{m(x|\pi)}\right)-\ln\left(\frac{\pi_0(\theta)}{q(\theta)}\right)\right). \label{derivative1}
\end{eqnarray}
As
$$\frac{{\partial}}{\partial{\epsilon}}\ln \left(\frac{c(\epsilon)}{m(x|\pi)}\right)=E_{\pi_0(\theta|x)}\left[\ln\left(\frac{\pi_0(\theta)}{q(\theta)}\right)\right]$$
(Dey and Birmiwal, 1994, Theorem 3.2), we get
\begin{eqnarray*}
\frac{\partial \gamma}{\partial \epsilon}=a \gamma\left(E_{\pi_0(\theta|x)}\left[\ln\left(\frac{\pi_0(\theta)}{q(\theta)}\right)\right]-\ln\left(\frac{\pi_0(\theta)}{q(\theta)}\right)\right).
\end{eqnarray*}
Since $\gamma\bigg|_{\epsilon=0}=\pi_0(\theta|x)$, by    (\ref{2nd}) and (\ref{derivative1}), it follows that $\int\frac{\partial \gamma}{\partial\epsilon}d\theta\bigg|_{\epsilon=0}=0$
and
\begin{eqnarray*}
\frac{{\partial^2}d}{\partial{\epsilon^2}}\bigg|_{\epsilon=0}&=&\int\frac{\partial^2\gamma}{\partial\epsilon^2}d\theta\bigg|_{\epsilon=0}.
\end{eqnarray*}
Now, by (\ref{derivative1}),
\begin{eqnarray*}
\frac{\partial^2 \gamma}{\partial \epsilon^2}&=&\frac{\partial}{\partial \epsilon}\left(a \gamma\left(E_{\pi_0(\theta|x)}\left[\ln\left(\frac{\pi_0(\theta)}{q(\theta)}\right)\right]-\ln\left(\frac{\pi_0(\theta)}{q(\theta)}\right)\right)\right)\\
&=&a \gamma\left(E_{\pi_0(\theta|x)}\left[\ln\left(\frac{\pi_0(\theta)}{q(\theta)}\right)\right]-\ln\left(\frac{\pi_0(\theta)}{q(\theta)}\right)\right)^2.
\end{eqnarray*}
Using the $\gamma\bigg|_{\epsilon=0}=\pi_0(\theta|x)$ one more time, we obtain
\begin{equation*}
    \frac{{\partial^2}d}{\partial{\epsilon}^2}\bigg|_{\epsilon=0}=\int\frac{{\partial}^2{\gamma}}{\partial{\epsilon}^2}\bigg|_{\epsilon=0}d\theta=aVar_{{\pi_0}(\theta|x)}\left[\ln\left(\frac{q(\theta)}{{\pi_0}(\theta)}\right)\right].
\end{equation*}

\endproof

\section{Examples} In this section, the derived results are explained through three examples: the Bernoulli model, the multinomial model and the location normal model. In each example, the curvature values for the two classes \eqref{contaminated} and \eqref{geometric} are reported. Additionally, in Example 1, we computed R\'enyi divergence  between $\pi(\theta|x)$ and $\pi_{0}(\theta|x)$ and reported the calibrated value $p$ as described in (\ref{interpret1}) and (\ref{interpret2}). 
Recall that,  curvature values close to zero  indicate robustness of the used prior whereas larger values suggest lack of robustness. On the other hand, values of  $p$ close to 0.5 suggest robustness whereas   values of $p$ close to 1  means absence of robustness.

\medskip
\noindent {\textbf{Example 1 (Bernoulli Model).}}  \label{example1}
Suppose $x=(x_1, \ldots, x_n)$ is a sample from a Bernoulli distribution with a parameter $\theta$. Let the prior ${\pi_0}(\theta)$ be $B$eta$(\alpha,\beta)$. That is, $$\pi_0(\theta)=\frac{\Gamma(\alpha+\beta)}{\Gamma(\alpha)\Gamma(\beta)} \theta^{\alpha-1}(1-\theta)^{\beta-1}.$$ Thus, ${\pi_0}(\theta|x_1, \ldots, x_n)$ is
\begin{eqnarray}
B\mbox{eta}\left(\alpha+t,\beta+n-t\right),\label{post_example1}
\end{eqnarray}
where $t=\sum_{i=1}^n x_i.$ Let $q(\theta)$ be $B$eta$(c\alpha,c\beta)$ for $c>0$.

Now consider  the  sample $x=( 0, 0, 1, 1, 0, 1, 1, 1, 1,0, 0, 0, 1, 0, 1, 0, 1, 1, 0, 1)$ of size $n=20$ generated from $B$ernoulli$(0.5)$. For comparison purposes, we consider several values of $\alpha, \beta$ and $c$. Although it is possible to find exact formulas of the curvature by some algebraic manipulation, it looks more convenient to use a Monte Carlo approach in this example. First, we sample $\theta^{(s)}, s=1,\ldots, 10^6$, from the posterior distribution (\ref{post_example1}). Then we compute the variance of $q(\theta^{(s)})/\pi_0(\theta^{(s)})$ and the variance of $\ln\left(q(\theta^{(s)})/\pi_0(\theta^{(s)})\right)$. This can be implemented straightforwardly in \textbf{\textsf{R}}. The values of  the curvature for both classes \eqref{contaminated} and \eqref{geometric} are reported in Table \ref{tab1}. Remarkably, for the cases when $\alpha=\beta=1$ (uniform prior on $[0,1]$) and $\alpha=\beta=0.5$ (Jeffreys' prior), the curvature values are prominently small.

\begin{table}[htbp]
  \centering
  \setlength{\tabcolsep}{4.5 mm}
      \caption{Values of the local curvature for two classes $\Gamma_{a}$ and $\Gamma_{g}$ for a sample generated from Bernoulli(0.5).} \label{tab1}
      \scalebox{0.76}{
    \begin{tabular}[c]{llllllll}
    \toprule
\multirow{2}[3]{*}{$\begin{pmatrix}
\alpha\\ \beta
\end{pmatrix}$}  &\multirow{2}[3]{*}{$c$}&\multicolumn{2}{c}{$a=0.5$}&\multicolumn{2}{c}{$a=1$}&\multicolumn{2}{c}{$a=2$}
\\\cmidrule(lr){3-4}\cmidrule(lr){5-6}\cmidrule(lr){7-8}
&&$C_{a}^{\Gamma_{a}}$&$C_{a}^{\Gamma_{g}}$&$C_{a}^{\Gamma_{a}}$&$C_{a}^{\Gamma_{g}}$&$C_{a}^{\Gamma_{a}}$&$C_{a}^{\Gamma_{g}}$
\\\hline
\multirow{2}[3]{*}{$\begin{pmatrix}
0.5\\ 0.5
\end{pmatrix}$}&0.5&$8\times10^{-5}$&$0.0002$&$0.0001$&$0.0004$&$0.0003$&$0.0008$\\
&1&$0$&$0$&$0$&$0$&$0$&$0$\\
&1.5&$0.0003$&$0.0002$&$0.0006$&$0.0004$&$0.0013$&$0.0008$\\
&3&$0.0098$&$0.0033$&$0.0196$&$0.0067$&$0.0393$&$0.0135$\\
&5&$0.0531$&$0.0135$&$0.1062$&$0.0271$&$0.2125$&$0.0543$\\\hline
\multirow{2}[3]{*}{$\begin{pmatrix}
1\\ 1
\end{pmatrix}$}&0.5&$0.0003$&$0.0007$&$0.0007$&$0.0015$&$0.0014$&$0.0030$\\
&1&$0$&$0$&$0$&$0$&$0$&$0$\\
&1.5&$0.0010$&$0.0007$&$0.0021$&$0.0015$&$0.0042$&$0.0030$\\
&3&$0.0241$&$0.0121$&$0.0483$&$0.0243$&$0.0967$&$0.0486$\\
&5&$0.1065$&$0.0486$&$0.2130$&$0.0972$&$0.4260$&$0.1945$\\\hline
\multirow{2}[3]{*}{$\begin{pmatrix}
1\\ 3
\end{pmatrix}$}&0.5&$0.0265$&$0.0235$&$0.0530$&$0.0470$&$0.1060$&$0.0941$\\
&1&$0$&$0$&$0$&$0$&$0$&$0$\\
&1.5&$0.0171$&$0.0235$&$0.0342$&$0.0470$&$0.0684$&$0.0941$\\
&3&$0.1061$&$0.3767$&$0.2122$&$0.7535$&$0.4244$&$1.5070$\\
&5&$0.1660$&$1.5070$&$0.3320$&$3.0141$&$0.6641$&$6.0282$\\\hline
\multirow{2}[3]{*}{$\begin{pmatrix}
3\\ 1
\end{pmatrix}$}&0.5&$0.0089$&$0.0113$&$0.0179$&$0.0227$&$0.0133$&$0.0454$\\
&1&$0$&$0$&$0$&$0$&$0$&$0$\\
&1.5&$0.0108$&$0.0113$&$0.0216$&$0.0227$&$0.0433$&$0.0454$\\
&3&$0.1162$&$0.1819$&$0.2324$&$0.3638$&$0.4648$&$0.7277$\\
&5&$0.2774$&$0.7277$&$0.5548$&$1.4555$&$1.1096$&$2.9110$\\
    \bottomrule
     \end{tabular}
  } %
\end{table}%

While it is  easier to quantify  the curvature based on Theorem \ref{theorem1} and Theorem  \ref{theorem2}, in this example, for comparison purposes, we computed  R\'enyi divergence between $\pi(\theta|x)$ and $\pi_0(\theta|x)$ under class \eqref{contaminated} and class \eqref{geometric}. It can be shown that, under class \eqref{contaminated} in \eqref{posterior1},
$\pi(\theta|x)=\lambda(x)B\mbox{eta}\left(\alpha+t,\beta+n-t\right)+(1-\lambda(x))B\mbox{eta}\left(c\alpha+t,c\beta+n-t\right),$
where
\begin{align*}
\lambda(x)=\frac{(1-\epsilon)\frac{\Gamma(\alpha+\beta)}{\Gamma(\alpha)\Gamma(\beta)}\frac{\Gamma(\alpha+t)\Gamma(\beta-t+n)}{\Gamma(\alpha+\beta+n)}}{(1-\epsilon)\frac{\Gamma(\alpha+\beta)}{\Gamma(\alpha)\Gamma(\beta)}\frac{\Gamma(\alpha+t)\Gamma(\beta-t+n)}{\Gamma(\alpha+\beta+n)}+\epsilon\frac{\Gamma(c\alpha+c\beta)}{\Gamma(c\alpha)\Gamma(c\beta)}\frac{\Gamma(c\alpha+t)\Gamma(c\beta-t+n)}{\Gamma(c\alpha+c\beta+n)}}.
\end{align*}
Also, from \eqref{posterior2}, it can be easily concluded that the posterior $\pi(\theta|x)$   under class \eqref{geometric} is obtained as
\begin{align*}
\pi(\theta|x)&=K \times \frac{\theta^{t}(1-\theta)^{n-t}\left[B\mbox{eta}\left(\alpha,\beta\right)\right]^{1-\epsilon}\left[B\mbox{eta}\left(c\alpha,c\beta\right)\right]^{\epsilon}}{\left[ \frac{\Gamma(\alpha+\beta)}{\Gamma(\alpha)\Gamma(\beta)} \right]^{(1-\epsilon)}\left[
\frac{\Gamma(c\alpha+c\beta)}{\Gamma(c\alpha)\Gamma(c\beta)}
 \right]^{\epsilon}},
 \end{align*}
$$K=\frac{\Gamma(t+(1-\epsilon)(\alpha-1)+\epsilon(c\alpha-1)+1)\Gamma(n-t+(1-\epsilon)
(\beta-1)+\epsilon(c\beta-1)+1)}{\Gamma((1-\epsilon)(\alpha+\beta-2)+\epsilon(c\alpha+c\beta-2)+n+2)}.$$
Note that, since  $d(\pi(\theta|x),\pi_{0}(\theta|x))= \frac{1}{a-1} \ln\left(E_{\pi_{0}(\theta|x)}\left[\left(\frac{\pi(\theta|x)}{\pi_{0}(\theta|x)}\right)^a\right]\right)$,  it possible  to  compute  the distance  based on a Monte Carlo approach. When $a=1$, $d(\pi(\theta|x),\pi_{0}(\theta|x))=E_{\pi_{0}(\theta|x)}\left[\frac{\pi(\theta|x)}{\pi_{0}(\theta|x)}\ln\left(\frac{\pi(\theta|x)}{\pi_{0}(\theta|x)}\right)\right]$, the Kullback-Leibler divergence.
We also  calibrated  R\'enyi divergence as described in (\ref{interpret1}) and (\ref{interpret2}). 
The results based on class \eqref{contaminated} and \eqref{geometric} are reported, respectively, in Table \ref{exm1under1} and Table \ref{exm1under2}.
\smallskip
\begin{table}[htbp]
  \centering
  \setlength{\extrarowheight}{-1mm}
  \setlength{\tabcolsep}{2.5 mm}
      \caption{Values of $d_0$ and $p$ in \eqref{interpret1} (for $a\neq 1$) and \eqref{interpret2} (for $a=1$) under class \eqref{contaminated} for a sample generated from Bernoulli(0.5).} \label{exm1under1}
      \scalebox{0.7}{
    \begin{tabular}[c]{llllllllllll}
    \toprule
\multirow{2}[3]{*}{$\begin{pmatrix}
\alpha\\ \beta
\end{pmatrix}$}  &\multirow{2}[3]{*}{$c$}&&\multicolumn{3}{c}{$a=0.5$}&\multicolumn{3}{c}{$a=1$}&\multicolumn{3}{c}{$a=2$}
\\\cmidrule(lr){4-6}\cmidrule(lr){7-9}\cmidrule(lr){10-12}
&&&$\epsilon=0.05$&$\epsilon=0.5$&$\epsilon=1$&$\epsilon=0.05$&$\epsilon=0.5$&$\epsilon=1$&$\epsilon=0.05$&$\epsilon=0.5$&$\epsilon=1$
\\\hline
\multirow{2}[3]{*}{$\begin{pmatrix}
0.5\\ 0.5
\end{pmatrix}$}&0.5&$d_0$&$2\times10^{-7}$&$4\times10^{-6}$&$9\times10^{-5}$&$5\times10^{-7}$&$3\times10^{-5}$&$0.0002$&$10^{-6}$&$7\times10^{-7}$&$0.0004$\\
&&$p$&(0.5003)&(0.5022)&(0.51)&(0.5005)&(0.5042)&(0.5107)&(0.5003)&(0.5041)&(0.5106)\\\cmidrule(lr){3-12}
&1&$d_0$&0&0&0&0&0&0&0&0&0\\
&&$p$&(0.5)&(0.5)&(0.5)&(0.5)&(0.5)&(0.5)&(0.5)&(0.5)&(0.5)\\\cmidrule(lr){3-12}
&1.5&$d_0$&$2\times10^{-6}$&$4\times10^{-5}$&$0.0001$&$2\times10^{-7}$&$5\times10^{-5}$&$0.0001$&$3\times10^{-7}$&$0.0001$&$0.0003$\\
&&$p$&(0.5013)&(0.5068)&(0.5104)&(0.5003)&(0.5054)&(0.5098)&(0.5003)&(0.5053)&(0.5096)\\\cmidrule(lr){3-12}
&3&$d_0$&$4\times10^{-6}$&$0.0004$&$0.0015$&$10^{-5}$&$0.0012$&$0.0028$&$3\times10^{-5}$&$0.0023$&$0.0054$\\
&&$p$&(0.5022)&(0.5204)&(0.5393)&(0.5031)&(0.5244)&(0.5379)&(0.5030)&(0.5239)&(0.5367)\\\cmidrule(lr){3-12}
&5&$d_0$&$5\times10^{-5}$&$0.0019$&$0.0055$&$0.0001$&$0.0048$&$0.0102$&$0.0002$&$0.0090$&$0.0181$\\
&&$p$&(0.5071)&(0.5437)&(0.5741)&(0.5074)&(0.5493)&(0.5711)&(0.5074)&(0.5476)&(0.5676)\\\hline
\multirow{2}[3]{*}{$\begin{pmatrix}
1\\ 1
\end{pmatrix}$}&0.5&$d_0$&$7\times10^{-7}$&$5\times10^{-5}$&$0.0003$&$10^{-6}$&$0.0001$&$0.0008$&$3\times10^{-6}$&$0.0002$&$0.0017$\\
&&$p$&(0.5007)&(0.5071)&(0.5193)&(0.5009)&(0.5083)&(0.5204)&(0.5007)&(0.5084)&(0.5207)\\\cmidrule(lr){3-12}
&1&$d_0$&0&0&0&0&0&0&0&0&0\\
&&$p$&(0.5)&(0.5)&(0.5)&(0.5)&(0.5)&(0.5)&(0.5)&(0.5)&(0.5)\\\cmidrule(lr){3-12}
&1.5&$d_0$&$2\times10^{-7}$&$7\times10^{-5}$&$0.0003$&$10^{-6}$&$0.0002$&$0.0006$&$2\times10^{-6}$&$0.0003$&$0.0013$\\
&&$p$&(0.5003)&(0.5084)&(0.5193)&(0.5008)&(0.5100)&(0.5185)&(0.5007)&(0.51)&(0.5180)\\\cmidrule(lr){3-12}
&3&$d_0$&$10^{-5}$&$0.0013$&$0.0050$&$5\times10^{-5}$&$0.0034$&$0.0092$&$0.0001$&$0.0065$&$0.0165$\\
&&$p$&(0.5042)&(0.5364)&(0.5706)&(0.5050)&(0.5416)&(0.5677)&(0.505)&(0.5405)&(0.5645)\\\cmidrule(lr){3-12}
&5&$d_0$&$8\times10^{-5}$&$0.0050$&$0.0167$&0.0002&0.0124&0.0297&0.0004&0.0225&0.0494\\
&&$p$&(0.5092)&(0.5708)&(0.6279)&(0.5107)&(0.5785)&(0.6201)&(0.5106)&(0.5755)&(0.6125)\\\hline
\multirow{2}[3]{*}{$\begin{pmatrix}
1\\ 3
\end{pmatrix}$}&0.5&$d_0$&$2\times10^{-5}$&$0.0032$&$0.0133$&$7\times10^{-5}$&$0.0067$&$0.0282$&$0.0001$&$0.0145$&$0.0623$\\
&&$p$&(0.5053)&(0.5565)&(0.6143)&(0.5059)&(0.5580)&(0.6171)&(0.5060)&(0.5604)&(0.6268)\\\cmidrule(lr){3-12}
&1&$d_0$&0&0&0&0&0&0&0&0&0\\
&&$p$&(0.5)&(0.5)&(0.5)&(0.5)&(0.5)&(0.5)&(0.5)&(0.5)&(0.5)\\\cmidrule(lr){3-12}
&1.5&$d_0$&$2\times10^{-5}$&0.0023&0.0104&$3\times10^{-5}$&0.0045&0.0199&$7\times10^{-5}$&0.0088&0.0370\\
&&$p$&(0.505)&(0.5484)&(0.6015)&(0.5044)&(0.5476)&(0.5989)&(0.5044)&(0.5472)&(0.5971)\\\cmidrule(lr){3-12}
&&$p$&(0.5081)&(0.5846)&(0.6878)&(0.5077)&(0.5834)&(0.6795)&(0.5077)&(0.5833)&(0.6793)\\\cmidrule(lr){3-12}
&3&$d_0$&0.0001&0.0175&01213&0.0002&0.0349&0.2125&0.0005&0.0691&0.3421\\
&&$p$&(0.5119)&(0.6308)&(0.8181)&(0.5115)&(0.6299)&(0.7942)&(0.5117)&(0.6337)&(0.8193)\\\cmidrule(lr){3-12}
&5&$d_0$&0.0002&0.0308&0.3423&0.0004&0.0638&0.5519&0.0008&0.1337&0.6003\\
&&$p$&(0.5145)&(0.6715)&(0.9536)&(0.5146)&(0.6731)&(0.9087)&(0.5144)&(0.6891)&(0.9535)\\\hline
\multirow{2}[3]{*}{$\begin{pmatrix}
3\\ 1
\end{pmatrix}$}&0.5&$d_0$&$7\times10^{-6}$&0.0012&0.0063&$2\times10^{-5}$&0.0027&0.0135&$5\times10^{-5}$&0.0057&0.0295\\
&&$p$&(0.5026)&(0.5356)&(0.5791)&(0.5036)&(0.5369)&(0.5816)&(0.5034)&(0.5379)&(0.5866)\\\cmidrule(lr){3-12}
&1&$d_0$&0&0&0&0&0&0&0&0&0\\
&&$p$&(0.5)&(0.5)&(0.5)&(0.5)&(0.5)&(0.5)&(0.5)&(0.5)&(0.5)\\\cmidrule(lr){3-12}
&1.5&$d_0$&$10^{-5}$&0.0013&0.0051&$2\times10^{-5}$&0.0025&0.0096&$4\times10^{-5}$&0.0048&0.0180\\
&&$p$&(0.5040)&(0.5364)&(0.5713)&(0.5034)&(0.5354)&(0.5692)&(0.5032)&(0.535)&(0.5674)\\\cmidrule(lr){3-12}
&3&$d_0$&0.0001&0.0139&0.0600&0.0002&0.0286&0.1054&0.0005&0.0505&0.1711\\
&&$p$&(0.5125)&(0.6168)&(0.7342)&(0.5117)&(0.6143)&(0.7180)&(0.5119)&(0.6137)&(0.7160)\\\cmidrule(lr){3-12}
&5&$d_0$&0.0003&0.0340&0.1724&0.0006&0.0657&0.2786&0.0012&0.1231&0.4062\\
&&$p$&(0.5196)&(0.68)&(0.865)&(0.5183)&(0.6754)&(0.8268)&(0.5177)&(0.6809)&(0.8539)\\
    \bottomrule
     \end{tabular}
  } %
\end{table}%
\begin{table}[htbp]
  \centering
  \setlength{\extrarowheight}{-1mm}
  \setlength{\tabcolsep}{2.5 mm}
      \caption{Values of $d_0$ and $p$ in \eqref{interpret1} (for $a\neq 1$) and \eqref{interpret2} (for $a=1$) under class \eqref{geometric} for a sample generated from Bernoulli(0.5).} \label{exm1under2}
      \scalebox{0.7}{
    \begin{tabular}[c]{llllllllllll}
    \toprule
\multirow{2}[3]{*}{$\begin{pmatrix}
\alpha\\ \beta
\end{pmatrix}$}  &\multirow{2}[3]{*}{$c$}&&\multicolumn{3}{c}{$a=0.5$}&\multicolumn{3}{c}{$a=1$}&\multicolumn{3}{c}{$a=2$}
\\\cmidrule(lr){4-6}\cmidrule(lr){7-9}\cmidrule(lr){10-12}
&&&$\epsilon=0.05$&$\epsilon=0.5$&$\epsilon=1$&$\epsilon=0.05$&$\epsilon=0.5$&$\epsilon=1$&$\epsilon=0.05$&$\epsilon=0.5$&$\epsilon=1$
\\\hline
\multirow{2}[3]{*}{$\begin{pmatrix}
0.5\\ 0.5
\end{pmatrix}$}&0.5&$d_0$&$4\times10^{-7}$&$2\times10^{-5}$&$9\times10^{-5}$&$10^{-6}$&$5\times10^{-5}$&0.0002&$2\times10^{-6}$&0.0001&0.0004\\
&&$p$&(0.5007)&(0.5043)&(0.51)&(0.5007)&(0.5054)&(0.5107)&(0.5007)&(0.5053)&(0.5106)\\\cmidrule(lr){3-12}
&1&$d_0$&0&0&0&0&0&0&0&0&0\\
&&$p$&(0.5)&(0.5)&(0.5)&(0.5)&(0.5)&(0.5)&(0.5)&(0.5)&(0.5)\\\cmidrule(lr){3-12}
&1.5&$d_0$&$2\times10^{-6}$&$2\times10^{-5}$&0.0001&$3\times10^{-8}$&$4\times10^{-5}$&0.0001&$6\times10^{-8}$&$9\times10^{-5}$&0.0003\\
&&$p$&(0.5014)&(0.5053)&(0.5103)&(0.5001)&(0.5048)&(0.5098)&(0.5)&(0.505)&(0.5096)\\\cmidrule(lr){3-12}
&3&$d_0$&$10^{-6}$&0.0004&0.0015&$6\times10^{-6}$&0.0007&0.0028&$10^{-5}$&0.0014&0.0054\\
&&$p$&(0.5013)&(00.5203)&(0.5390)&(0.5017)&(0.5195)&(0.5379)&(0.5014)&(0.5191)&(0.5367)\\\cmidrule(lr){3-12}
&5&$d_0$&$9\times10^{-6}$&0.0015&0.0055&$2\times10^{-5}$&0.0028&0.0102&$5\times10^{-5}$&0.0054&0.0181\\
&&$p$&(0.5030)&(0.5390)&(0.5738)&(0.5038)&(0.5379)&(0.5711)&(0.5036)&(0.5367)&(0.5676)\\\hline
\multirow{2}[3]{*}{$\begin{pmatrix}
1\\ 1
\end{pmatrix}$}&0.5&$d_0$&$10^{-6}$&$7\times10^{-7}$&0.0003&$2\times10^{-6}$&0.0002&0.0008&$5\times10^{-6}$&0.0004&0.0017\\
&&$p$&(0.5011)&(0.5087)&(0.5193)&(0.5012)&(0.5101)&(0.5204)&(0.5011)&(0.5103)&(0.5207)\\\cmidrule(lr){3-12}
&1&$d_0$&0&0&0&0&0&0&0&0&0\\
&&$p$&(0.5)&(0.5)&(0.5)&(0.5)&(0.5)&(0.5)&(0.5)&(0.5)&(0.5)\\\cmidrule(lr){3-12}
&1.5&$d_0$&$6\times10^{-8}$&$5\times10^{-5}$&0.0003&$8\times10^{-7}$&0.0001&0.0006&$10^{-6}$&0.0003&0.0013\\
&&$p$&(0.5)&(0.5077)&(0.5193)&(0.5006)&(0.5093)&(0.5185)&(0.5007)&(0.5093)&(0.5180)\\\cmidrule(lr){3-12}
&3&$d_0$&$8\times10^{-6}$&0.0009&0.0050&$2\times10^{-5}$&0.0026&0.0092&$5\times10^{-5}$&0.0048&0.0165\\
&&$p$&(0.5027)&(0.5309)&(0.5706)&(0.5035)&(0.5360)&(0.5677)&(0.5037)&(0.535)&(0.5645)\\\cmidrule(lr){3-12}
&5&$d_0$&$3\times10^{-5}$&0.0035&0.0167&0.0001&0.0092&0.0297&0.0002&0.0165&0.0494\\
&&$p$&(0.5062)&(0.5596)&(0.6279)&(0.5074)&(0.5677)&(0.6201)&(0.5073)&(0.5645)&(0.6125)\\\hline
\multirow{2}[3]{*}{$\begin{pmatrix}
1\\ 3
\end{pmatrix}$}&0.5&$d_0$&$2\times10^{-5}$&0.0030&0.0133&$6\times10^{-5}$&0.0064&0.0282&0.0001&0.0135&0.0623\\
&&$p$&(0.505)&(0.5555)&(0.6143)&(0.5056)&(0.5566)&(0.6171)&(0.5054)&(0.5583)&(0.6268)\\\cmidrule(lr){3-12}
&1&$d_0$&0&0&0&0&0&0&0&0&0\\
&&$p$&(0.5)&(0.5)&(0.5)&(0.5)&(0.5)&(0.5)&(0.5)&(0.5)&(0.5)\\\cmidrule(lr){3-12}
&1.5&$d_0$&$3\times10^{-5}$&0.0028&0.0104&$5\times10^{-5}$&0.0053&0.0199&0.0001&0.0103&0.0370\\
&&$p$&(0.5059)&(0.5527)&(0.6015)&(0.5022)&(0.5517)&(0.5989)&(0.5053)&(0.5509)&(0.5971)\\\cmidrule(lr){3-12}
&3&$d_0$&0.0004&0.0373&0.1213&0.0008&0.0690&0.2125&0.0017&0.1210&0.3421\\
&&$p$&(0.5216)&(0.6878)&(0.8181)&(0.5211)&(0.6795)&(0.7942)&(0.5209)&(0.6793)&(0.8193)\\\cmidrule(lr){3-12}
&5&$d_0$&0.0018&0.1213&0.3423&0.0034&0.2125&0.5519&0.0067&0.3421&0.6003\\
&&$p$&(0.5425)&(0.8181)&(09536)&(0.5417)&(0.7942)&(0.9087)&(0.5411)&(0.8193)&(0.9535)\\\hline
\multirow{2}[3]{*}{$\begin{pmatrix}
3\\ 1
\end{pmatrix}$}&0.5&$d_0$&$10^{-5}$&0.0014&0.0063&$3\times10^{-5}$&0.0031&0.0135&$6\times10^{-5}$&0.0065&0.0295\\
&&$p$&(0.5031)&(0.5381)&(0.5791)&(0.5040)&(0.5394)&(0.5816)&(0.5039)&(0.5403)&(0.5866)\\\cmidrule(lr){3-12}
&1&$d_0$&0&0&0&0&0&0&0&0&0\\
&&$p$&(0.5)&(0.5)&(0.5)&(0.5)&(0.5)&(0.5)&(0.5)&(0.5)&(0.5)\\\cmidrule(lr){3-12}
&1.5&$d_0$&$10^{-5}$&0.0014&0.0052&$2\times10^{-5}$&0.0025&0.0096&$4\times10^{-5}$&0.0049&0.0180\\
&&$p$&(0.5041)&(0.5376)&(0.5720)&(0.5034)&(0.5359)&(0.5692)&(0.5033)&(0.5353)&(0.5674)\\\cmidrule(lr){3-12}
&3&$d_0$&0.0002&0.0185&0.0604&0.0004&0.0338&0.1054&0.0008&0.0596&0.1711\\
&&$p$&(0.5153)&(0.6341)&(0.735)&(0.5145)&(0.6278)&(0.7180)&(0.5143)&(0.6239)&(0.7160)\\\cmidrule(lr){3-12}
&5&$d_0$&0.0008&0.0604&0.1724&0.0016&0.1054&0.2786&0.0032&0.1711&0.4074\\
&&$p$&(0.53)&(0.735)&(0.865)&(0.5289)&(0.7180)&(0.8268)&(0.5284)&(0.7160)&(0.8545)\\
    \bottomrule
     \end{tabular}
  } %
\end{table}%

Note that, from  (\ref{fisher}), by multiplying the curvature value in  Table \ref{tab1} by $\epsilon^2/2$, one may get the value of the corresponding distance in Table \ref{exm1under1} and Table \ref{exm1under2}.  For instance, setting $\alpha=1, \beta =1, c=0.5, a=0.5$ in Table  \ref{tab1}, gives $C_{a}^{\Gamma_{a}}=0.0265$. The corresponding  distance is $0.0265 \times 0.5^2/2= 0.0033$, which close to the one reported in Table  \ref{exm1under1}.

Now we consider the Australian AIDS survival data, available in the \textbf{\textsf{R}} package ``\textsf{Mass}". There are  2843 patients diagnosed with AIDS in Australia before 1 July 1991. The data frame contains the following columns: state, sex, date of diagnosis,  date of death at end of observation, status (``$A$" (alive) or ``$D$" (dead) at end of observation), reported transmission category, and age at diagnosis.  Now, we consider  the values of column status,  then, under prior distribution given above,  the values of  the curvatures for two classes \eqref{contaminated} and \eqref{geometric} are summarized in Table \ref{tab-real}.
\begin{table}[htbp]
  \centering
  \setlength{\tabcolsep}{4 mm}
      \caption{Values of the local curvature for the two classes $\Gamma_{a}$ and $\Gamma_{g}$ for the real data set AIDS.} \label{tab-real}
      \scalebox{0.76}{
    \begin{tabular}[c]{llllllll}
    \toprule
\multirow{2}[3]{*}{$\begin{pmatrix}
\alpha\\ \beta
\end{pmatrix}$}  &\multirow{2}[3]{*}{$c$}&\multicolumn{2}{c}{$a=0.5$}&\multicolumn{2}{c}{$a=1$}&\multicolumn{2}{c}{$a=2$}
\\\cmidrule(lr){3-4}\cmidrule(lr){5-6}\cmidrule(lr){7-8}
&&$C_{a}^{\Gamma_{a}}$&$C_{a}^{\Gamma_{g}}$&$C_{a}^{\Gamma_{a}}$&$C_{a}^{\Gamma_{g}}$&$C_{a}^{\Gamma_{a}}$&$C_{a}^{\Gamma_{g}}$
\\\hline
\multirow{2}[3]{*}{$\begin{pmatrix}
0.5\\ 0.5
\end{pmatrix}$}&0.5&$9\times10^{-7}$&$2\times10^{-6}$&$10^{-6}$&$5\times10^{-6}$&$3\times10^{-6}$&$10^{-5}$\\
&1&0&0&0&0&0&0\\
&1.5&$4\times10^{-6}$&$2\times10^{-6}$&$8\times10^{-6}$&$5\times10^{-6}$&$10^{-5}$&$10^{-5}$\\
&3&0.0001&$4\times10^{-5}$&0.0003&$8\times10^{-5}$&0.0006&0.0001\\
&5&0.0009&0.0001&0.0019&0.0003&0.0038&0.0006\\
\hline
\multirow{2}[3]{*}{$\begin{pmatrix}
1\\ 1
\end{pmatrix}$}&0.5&$4\times10^{-6}$&$10^{-5}$&$9\times10^{-6}$&$2\times10^{-5}$&$10^{-5}$&$4\times10^{-5}$\\
&1&0&0&0&0&0&0\\
&1.5&$10^{-5}$&$10^{-5}$&$3\times10^{-5}$&$2\times10^{-5}$&$6\times10^{-5}$&$4\times10^{-5}$\\
&3&0.0004&0.0001&0.0009&0.0003&0.0018&0.0006\\
&5&0.0025&0.0006&0.0051&0.0013&0.0102&0.0027\\
\hline
\multirow{2}[3]{*}{$\begin{pmatrix}
1\\ 3
\end{pmatrix}$}&0.5&0.0005&0.0004&0.0010&0.0008&0.0021&0.0016\\
&1&0&0&0&0&0&0\\
&1.5&0.0002&0.0004&0.0004&0.0008&0.0008&0.0016\\
&3&0.0002&0.0064&0.0004&0.0129&0.0009&0.0259\\
&5&$10^{-5}$&0.0259&$3\times10^{-5}$&0.0518&$7\times10^{-5}$&0.1037\\
\hline
\multirow{2}[3]{*}{$\begin{pmatrix}
3\\ 1
\end{pmatrix}$}&0.5&$2\times10^{-5}$&$5\times10^{-5}$&$5\times10^{-5}$&0.0001&0.0001&0.0002\\
&1&0&0&0&0&0&0\\
&1.5&$6\times10^{-5}$&$5\times10^{-5}$&0.0001&0.0001&0.0002&0.0002\\
&3&0.0014&0.0008&0.0029&0.0016&0.0058&0.0032\\
&5&0.0054&0.0032&0.0108&0.0064&0.0216&0.0129\\
    \bottomrule
     \end{tabular}
  }
\end{table}%

\textbf{Example 2 (Multinomial model).} \label{example2}
Suppose that $x=(x_1,x_2,\ldots, x_k)$ is an observation  from a multinomial distribution with parameters $(N,(\theta_1,\ldots,\theta_k))$, where $\sum_{i=1}^{k}x_i=N$ and $\sum_{i=1}^{k}\theta_i=1$.  Let the prior $\pi_0(\theta_1,\ldots,\theta_k)$  be $D$irichlet$(\alpha_1,$ $\ldots,\alpha_k)$. Then $\pi_0(\theta_1,\ldots,\theta_k|x)$ is $D\mbox{irichlet}(\alpha_1+x_1,\ldots,\alpha_k+x_k).$

Let $q(\theta_1,\ldots,\theta_k)\sim D\mbox{irichlet}(c\alpha_1,\ldots,c\alpha_k)$. We consider the  observation  $x=(6,4,5,5)$ generated  from $M$ultinomial$(20,(1/4,1/4,1/4,1/4))$. As in Example 1, we use Monte Carlo approach to compute curvature values. Table \ref{tab2} reports values of the curvature for different values of $\alpha_1,\ldots,\alpha_k$ and $c$. Clearly, when $c=1$, the curvature values are 0. Also, for the cases when $\alpha_1=\alpha_2=\alpha_3=\alpha_4=1$ (uniform prior over $[0,1]^4$) and $\alpha_1=\alpha_2=\alpha_3=\alpha_4=0.5$ (Jeffreys' prior), the curvature values are prominently small.

\begin{table}[htbp]
  \centering
  \setlength{\tabcolsep}{4.5 mm}
      \caption{Values of the local curvature for two classes $\Gamma_{a}$ and $\Gamma_{g}$ for a sample generated from Mn(20,(1/4,1/4,1/4,1/4)).}\label{tab2}
      \scalebox{0.76}{
    \begin{tabular}[c]{llllllll}
    \toprule
\multirow{2}[3]{*}{$\Biggl(\begin{smallmatrix}
\alpha_{1}\\ \vdots\\ \alpha_{4}
\end{smallmatrix}\Biggr)$}  &\multirow{2}[3]{*}{$c$}&\multicolumn{2}{c}{$a=0.5$}&\multicolumn{2}{c}{$a=1$}&\multicolumn{2}{c}{$a=2$}
\\\cmidrule(lr){3-4}\cmidrule(lr){5-6}\cmidrule(lr){7-8}
&&$C_{a}^{\Gamma_{a}}$&$C_{a}^{\Gamma_{g}}$&$C_{a}^{\Gamma_{a}}$&$C_{a}^{\Gamma_{g}}$&$C_{a}^{\Gamma_{a}}$&$C_{a}^{\Gamma_{g}}$
\\\hline
\multirow{4}[3]{*}{$\begin{pmatrix}
0.25\\ 0.25\\ 0.25\\ 0.25\\
\end{pmatrix}$}&0.5&$2\times10^{-5}$&$0.0006$&$5\times10^{-5}$&$0.0012$&$0.0001$&$0.0024$\\
&1&$0$&$0$&$0$&$0$&$0$&$0$\\
&1.5&$0.0031$&$0.0006$&$0.0062$&$0.0012$&$0.0124$&$0.0024$\\
&3&$0.5285$&$0.0097$&$1.0570$&$0.0195$&$2.1141$&$0.0390$\\
&5&$8.4050$&$0.0301$&$16.816$&$0.0780$&$33.632$&$0.1560$\\\hline
\multirow{4}[3]{*}{$\begin{pmatrix}
0.5\\ 0.5\\ 0.5\\ 0.5\\
\end{pmatrix}$}&0.5&$0.0001$&$0.0021$&$0.0003$&$0.0043$&$0.0004$&$0.0087$\\
&1&$0$&$0$&$0$&$0$&$0$&$0$\\
&1.5&$0.0080$&$0.0021$&$0.0161$&$0.0043$&$0.0323$&$0.0087$\\
&3&$0.7706$&$0.0349$&$1.5413$&$0.0699$&$3.0826$&$0.1398$\\
&5&$8.0246$&$0.1398$&$16.049$&$0.2797$&$32.098$&$0.5595$\\\hline
\multirow{4}[3]{*}{$\begin{pmatrix}
1\\ 1\\ 1\\ 1\\
\end{pmatrix}$}&0.5&$0.0008$&$0.0071$&$0.0017$&$0.0142$&$0.0035$&$0.0284$\\
&1&$0$&$0$&$0$&$0$&$0$&$0$\\
&1.5&$0.0185$&$0.0071$&$0.0370$&$0.0142$&$0.0741$&$0.0284$\\
&3&$0.9799$&$0.1137$&$1.9598$&$0.2274$&$3.9196$&$0.4549$\\
&5&$6.7661$&$0.4549$&$13.532$&$0.9098$&$27.064$&$1.8197$\\\hline
\multirow{4}[3]{*}{$\begin{pmatrix}
2\\ 1\\ 1\\ 1\\
\end{pmatrix}$}&0.5&$0.0018$&$0.0120$&$0.0037$&$0.0240$&$0.0074$&$0.0480$\\
&1&$0$&$0$&$0$&$0$&$0$&$0$\\
&1.5&$0.0270$&$0.0120$&$0.0540$&$0.0240$&$0.1081$&$0.0480$\\
&3&$1.1052$&$0.1923$&$2.2104$&$0.3847$&$4.4209$&$0.7695$\\
&5&$6.3984$&$0.7695$&$12.796$&$1.5390$&$25.593$&$3.0780$\\
    \bottomrule
     \end{tabular}
  }
\end{table}%

\smallskip

\textbf{Example 3 (Location normal model).} \label{example3}
Suppose that $x=(x_1,x_2,\ldots, x_n)$ is a sample  from $N(\theta,1)$ distribution with $\theta \in \mathbb{R}^1$.  Let the prior $\pi_{0}(\theta)$ of $\theta$ be $N(\theta_0,\sigma_0^2)$. Then
\begin{align}
\begin{split}
\label{normal}
\pi_{0}(\theta|x)\sim \mathcal{N}\left(\mu_x,\sigma^2_{x}\right),
\end{split}
\end{align}
\begin{eqnarray*}
\mu_x=\left(\frac{\theta_0}{\sigma_0^2}+n\bar{x}\right)\left(\frac{1}{\sigma_0^2}+n\right)^{-1}\ \ \text{and} \ \ \sigma^2_{x}=\left(\frac{1}{\sigma_0^2}+n\right)^{-1}. 
\end{eqnarray*}
Let $q(\theta)\sim \mathcal{N}(c\theta_0,\sigma_{0}^{2})$, $c>0$. Due to some interesting theoretical properties in this example, we  present the exact formulas of the curvature for class \eqref{contaminated} and class \eqref{geometric}. We have
\begin{align}
\begin{split}
 \nonumber   \frac{q(\theta)}{\pi_0(\theta)}=\exp\left\{\frac{\theta_0\theta(c-1)+0.5\theta_0^2(1-c^2)}{\sigma_0^2}\right\}.
\end{split}
\end{align}
Therefore, for the class (\ref{contaminated}), we have
\begin{eqnarray*}
Var_{{\pi_0}(\theta|x)}\left[\frac{q(\theta)}{\pi_0(\theta)}\right]&=&E_{{\pi_0}(\theta|x)}\left[\left(\frac{q(\theta)}{\pi_0(\theta)}\right)^2\right]
-\left(E_{{\pi_0}(\theta|x)}\left[\frac{q(\theta)}{\pi_0(\theta)}\right]\right)^2\\
&=&\exp\left\{\frac{\theta_0^2(1-c^2)}{\sigma_0^2}\right\}\bigg[M_{{\pi_0}(\theta|x)}\left(\frac{2\theta_0(c-1)}{\sigma_0^2}\right)-\\
    &&\left(M_{{\pi_0}(\theta|x)}\left(\frac{\theta_0(c-1)}{\sigma_0^2}\right)\right)^2\bigg],
\end{eqnarray*}
where $M_{{\pi_0}(\theta|x)}(t)$ is the moment generating function with respect to the density $\pi_{0}(\theta|x)$. Thus,
\begin{eqnarray*}
Var_{{\pi_0}(\theta|x)}\left[\frac{q(\theta)}{\pi_0(\theta)}\right]&=&
\exp\left\{\frac{\theta_0^2(1-c^2)}{\sigma_0^2}\right\}\bigg[\exp\bigg\{\frac{2\theta_0(c-1)\mu_x}{\sigma_0^2}+\\
&&\frac{2\theta_0^2(c-1)^2 \sigma_x^2}{\sigma_0^4}\bigg\}-\exp\bigg\{\frac{2\theta_0(c-1)\mu_x}{\sigma_0^2}+\\
&&\frac{\theta_{0}^2(c-1)^2\sigma_x^2}{\sigma_0^4}\bigg\}\bigg].
\end{eqnarray*}
On the other hand, for  the geometric contaminated class, we have
\begin{align*}
\ln\left(\frac{q(\theta)}{\pi_0(\theta)}\right)=\frac{\theta_0\theta(c-1)+0.5\theta_0^2(1-c^2)}{\sigma_0^2}.
\end{align*}
Thus, by (\ref{normal}), we get
\begin{eqnarray}
\nonumber Var_{{\pi_0}(\theta|x)}\left[\ln\left(\frac{q(\theta)}{\pi_0(\theta)}\right)\right]&=&
\frac{\theta_0^2(c-1)^2}{\sigma_0^4} Var_{{\pi_0}(\theta|x)}\left[\theta\right]\\
\nonumber &=&\frac{\theta_0^2(c-1)^2}{\sigma_0^4} \sigma^2_x\\
&=&\frac{\theta_0^2(c-1)^2}{\sigma_0^4} \left(\frac{1}{\sigma_0^2}+n\right)^{-1}. \label{example2_geometric}
\end{eqnarray}
Interestingly, from (\ref{example2_geometric}), $Var_{{\pi_0}(\theta|x)}\left[\ln\left(\frac{q(\theta)}{\pi_0(\theta)}\right)\right]$ depends on the sample only through its size $n$. As $n \to \infty$ or $\sigma_0 \to \infty$, $Var_{{\pi_0}(\theta|x)}\left[\ln\left(\frac{q(\theta)}{\pi_0(\theta)}\right)\right] \to 0,$ which indicates robustness. Also, as $\theta_{0} \to \infty$, $Var_{{\pi_0}(\theta|x)}\left[\ln\left(\frac{q(\theta)}{\pi_0(\theta)}\right)\right] \to \infty$ and no robustness will be found.

Now we consider a numerical example by generating a sample of size $n=20$ from $N(4,1)$ distribution. We obtain
\begin{quote}
$x=(3.37, 4.18, 3.16, 5.59, 4.32, 3.17, 4.48, 4.73,
  4.57, 3.69, 5.51, 4.38, 3.37,\newline
   1.78, 5.12, 3.95,
3.98, 4.94, 4.82, 4.59)$
\end{quote}
(with $t=\bar{x}=4.1905$). Table \ref{tab3} reports the values of the curvature for different values of $\theta_0, \sigma_0$ and $c$.

\begin{table}[htbp]
  \centering
  \setlength{\tabcolsep}{4.5 mm}
      \caption{Values of the local curvature for two classes $\Gamma_{a}$ and $\Gamma_{g}$ for a sample generated from N(4,1).} \label{tab3}
      \scalebox{0.76}{
    \begin{tabular}[c]{llllllll}
    \toprule
\multirow{2}[3]{*}{$\begin{pmatrix}
\theta_{0} \\ \sigma_{0}^{2}
\end{pmatrix}$}  &\multirow{2}[3]{*}{$c$}&\multicolumn{2}{c}{$a=0.5$}&\multicolumn{2}{c}{$a=1$}&\multicolumn{2}{c}{$a=2$}
\\\cmidrule(lr){3-4}\cmidrule(lr){5-6}\cmidrule(lr){7-8}
&&$C_{a}^{\Gamma_{a}}$&$C_{a}^{\Gamma_{g}}$&$C_{a}^{\Gamma_{a}}$&$C_{a}^{\Gamma_{g}}$&$C_{a}^{\Gamma_{a}}$&$C_{a}^{\Gamma_{g}}$
\\\hline
\multirow{2}[3]{*}{$\begin{pmatrix}
0.1\\ 0.1
\end{pmatrix}$}&0.5&$0.0001$&$0.0059$&$0.0002$&$0.0119$&$0.0004$&$0.0238$\\
&1&$0$&$0$&$0$&$0$&$0$&$0$\\
&1.5&$0.2908$&$0.0059$&$0.5816$&$0.0119$&$1.1633$&$0.0238$\\
&3&$498033.7$&$0.0953$&$996067.4$&$0.1907$&$1992135$&$0.3814$\\
&5&$8\times10^{12}$&$0.3814$&$10^{13}$&$0.7629$&$3\times10^{13}$&$1.5258$\\\hline
\multirow{2}[3]{*}{$\begin{pmatrix}
0.5\\ 1
\end{pmatrix}$}&0.5&$0.0002$&$0.0014$&$0.0004$&$0.0029$&$0.0009$&$0.0059$\\
&1&$0$&$0$&$0$&$0$&$0$&$0$\\
&1.5&$0.0081$&$0.0014$&$0.0162$&$0.0029$&$0.0325$&$0.0059$\\
&3&$10.629$&$0.0238$&$21.258$&$0.0476$&$42.517$&$0.0953$\\
&5&$2964.9$&$0.0935$&$2929.8$&$0.1907$&$11859.7$&$0.3814$\\\hline
\multirow{2}[3]{*}{$\begin{pmatrix}
0.5\\ 5
\end{pmatrix}$}&0.5&$4\times10^{-5}$&$5\times10^{-5}$&$8\times10^{-5}$&$0.0001$&$0.0001$&$0.0002$\\
&1&$0$&$0$&$0$&$0$&$0$&$0$\\
&1.5&$8\times10^{-5}$&$5\times10^{-5}$&$0.0001$&$0.0001$&$0.0003$&$0.0002$\\
&3&$0.0031$&$0.0009$&$0.0063$&$0.0019$&$0.0127$&$0.0038$\\
&5&$0.0288$&$0.0038$&$0.0576$&$0.0076$&$0.1152$&$0.0152$\\\hline
\multirow{2}[3]{*}{$\begin{pmatrix}
4\\ 5
\end{pmatrix}$}&0.5&$0.0001$&$0.0038$&$0.0029$&$0.0076$&$0.0059$&$0.0152$\\
&1&$0$&$0$&$0$&$0$&$0$&$0$\\
&1.5&$0.0020$&$0.0038$&$0.0040$&$0.0076$&$0.0080$&$0.0152$\\
&3&$3\times10^{-7}$&$0.0610$&$7\times10^{-7}$&$0.1220$&$10^{-6}$&$0.2441$\\
&5&$9\times10^{-23}$&$0.2441$&$10^{-22}$&$0.4882$&$3\times10^{-22}$&$0.9765$\\
    \bottomrule
     \end{tabular}
  }
\end{table}%

Clearly, for large values of $\sigma^2_0$, the value of the curvature is small, which is an indication of robustness. For instance, for $\mu_0=0.5$  in Table \ref{tab3}, that value of the curvature when $\sigma^2_0=5$ is much smaller than the value of the curvature when $\sigma^2_0=1$.

\smallskip

\section{Conclusions}  Measuring Bayesian robustness of two classes of contaminated priors is studied. The approach is based on computing the curvature of  R\'enyi divergence  between posterior distributions.  The method does not require specifying values for $\epsilon$ and its computation is straightforward. Examples illustrating the approach are considered.

%

\end{document}